\documentclass[12pt,reqno]{amsart}
\usepackage{fullpage}
\usepackage{mathpple}

\newtheorem{theorem}{Theorem}
\newtheorem{lemma}[theorem]{Lemma}

\newtheorem{cor}[theorem]{Corollary}

\renewcommand{\mod}[1]{{\ifmmode\text{\rm\ (mod~$#1$)}\else\discretionary{}{}{\hbox{ }}\rm(mod~$#1$)\fi}}
\newcommand{\ep}{\varepsilon}
\newsymbol\dnd 232D

\newcommand{\A}{{\mathcal A}}
\newcommand{\M}{{\mathcal M}}
\newcommand{\N}{{\mathcal N}}
\renewcommand{\P}{{\mathcal P}}

\begin{document}

\title{Simultaneous Inequalities among Values of the\\Euler $\phi$-function}
\author{Greg Martin} 
\address{Department of Mathematics \\ University of British Columbia \\ Room
121, 1984 Mathematics Road \\ Canada V6T 1Z2}
\email{gerg@math.ubc.ca}
\maketitle

\begin{abstract}
This paper concerns the values of the Euler $\phi$-function evaluated simultaneously on corresponding members of $k$ arithmetic progressions $a_1n+b_1$, $a_2n+b_2$, \dots, $a_kn+b_k$. Assuming the necessary condition that no two of the polynomials $a_ix+b_i$ are constant multiples of each other, we show that for any positive constant $C$, there are infinitely many integers $n$ for which simultaneously
\[
\frac{\phi(a_1n+b_1)}{\phi(a_2n+b_2)} > C, \quad
\frac{\phi(a_2n+b_2)}{\phi(a_3n+b_3)} > C, \quad\dots,\quad
\frac{\phi(a_{k-1}n+b_{k-1})}{\phi(a_kn+b_k)} > C.
\]
For example, there are infinitely many integers $n$ for which
\[
\phi(a_1n+b_1) > \phi(a_2n+b_2) > \dots > \phi(a_kn+b_k);
\]
in particular, there exist infinitely many strings of $k$ consecutive integers whose $\phi$-values are arranged from largest to smallest in any prescribed manner. Also, under the necessary condition $ad\ne bc$, any inequality of the form $\phi(an+b) > \phi(cn+d)$ infinitely often has $k$ consecutive solutions. In fact, we prove that the sets of solutions to these inequalities have positive lower density. All of these results hold as well when the Euler function $\phi$ is replaced by the sum-of-divisors function $\sigma$.
\end{abstract}

The Euler phi-function, like many multiplicative functions, is rather irregular because its values depend upon the factorizations of the corresponding arguments. Moreover, a value $\phi(n)$ can be smaller than the argument $n$ by an arbitrarily large factor. Thus when examining the values of $\phi(n)$ on multiplicatively unrelated arguments, such as values of different linear polynomials, we should not expect any one set of values to be consistently larger than another. One can construct situations where the sizes are heavily skewed in one direction; for example, the inequality $\phi(30n) < \phi(30n+1)$ holds for all $n\le 10^{1115}$ (see \cite{martin}). However, Newman \cite{newman} showed that any inequality of the form $\phi(an+b) < \phi(cn+d)$ holds infinitely often, provided only that the polynomials $an+b$ and $cn+d$ are not multiples of each other (an alternate proof was given in \cite{ABGU}).

We extend this result in several ways. First, we show that such inequalities can hold by an arbitrarily large factor; that is, we consider $\phi(an+b) < C\phi(cn+d)$ for some large constant $C>0$ rather than simply $\phi(an+b) < \phi(cn+d)$. Second, we extend the result to multiple simultaneous inequalities among many values of $\phi(n)$, as opposed to a single inequality between two values. Third, and perhaps most significantly, we quantify the number of solutions to such inequalities: rather than assert only that there are infinitely many solutions, we prove that the set of solutions actually has positive lower density.

Our main result is:

\begin{theorem}
Let $k\ge2$ be an integer, let $a_1,\dots,a_k$ be positive integers, and let $b_1,\dots,b_k$ be integers. Assume that $a_ib_j\ne a_jb_i$ for every $1\le i<j\le k$. Then for any positive constant $C$, the set of positive integers $n$ for which
\begin{equation}
\frac{\phi(a_1n+b_1)}{\phi(a_2n+b_2)} > C, \quad
\frac{\phi(a_2n+b_2)}{\phi(a_3n+b_3)} > C, \quad\dots,\quad
\frac{\phi(a_{k-1}n+b_{k-1})}{\phi(a_kn+b_k)} > C
\label{general.ineqs}
\end{equation}
has positive lower density.
In particular, the set of positive integers $n$ for which
\[
\phi(a_1n+b_1) > \phi(a_2n+b_2) > \dots > \phi(a_kn+b_k).
\]
has positive lower density.
\label{general.theorem}
\end{theorem}

We note that if $a_ib_j = a_jb_i$, then the ratio $\phi(a_in+b_i)/\phi(a_jn+b_j)$ can take only finitely many values. For example,  $\phi(3n+15)/\phi(n+5)$ equals either 3 or 2 according to whether or not 3 divides $n+5$, and so the inequality $\phi(3n+15) < \phi(n+5)$ has no solutions. Therefore the assumption that $a_ib_j\ne a_jb_i$ is necessary in Theorem \ref{general.theorem}.

We also remark that the hypotheses of Theorem \ref{general.theorem} are symmetric in the polynomials $a_ix+b_i$, and so the theorem implies that the set of positive integers $n$ satisfying any of the $k!$ prescribed orderings
\[
\phi(a_{\sigma(1)}n+b_{\sigma(1)}) > \phi(a_{\sigma(2)}n+b_{\sigma(2)}) > \dots > \phi(a_{\sigma(k)}n+b_{\sigma(k)})
\]
has positive lower density.

By taking $a_i=1$ and $(b_1,\dots,b_k)$ a permutation of $(1,\dots,k)$ in Theorem \ref{general.theorem}, we can immediately address any multiple inequalities of the form
\[
\newcommand{\gtlt}{\mathop{\genfrac{}{}{0pt}{1}{\textstyle>}{\textstyle<}}}
\phi(n+1) \gtlt \phi(n+2) \gtlt \cdots \gtlt \phi(n+k),
\]
where each inequality can be chosen to be $>$ or $<$ independently:

\begin{cor}
Let $k\ge2$ be an integer $\ep_1,\dots,\ep_{k-1}\in\{-1,1\}$. Then the set of positive integers $n$ for which
\begin{multline}
\ep_1(\phi(n+1)-\phi(n+2)) > 0,\quad \ep_2(\phi(n+2)-\phi(n+3)) > 0, \\
\dots,\quad \ep_{k-1}(\phi(n+k-1)-\phi(n+k)) > 0
\label{compare.Schinzel}
\end{multline}
has positive lower density.
\label{permutation.cor}
\end{cor}

For the sake of comparison, we note that Schinzel \cite{schinzel} proved that for every $(k-1)$-tuple of constants $\alpha_1,\alpha_2,\dots,\alpha_{k-1} \in [0,\infty]$, there exists a sequence $\{n_i\}$ such that
\[
\lim_{i\to\infty}\frac{\phi(n_i+1)}{\phi(n_i+2)} = \alpha_1,\quad
\lim_{i\to\infty}\frac{\phi(n_i+2)}{\phi(n_i+3)} = \alpha_2, \quad\dots,\quad
\lim_{i\to\infty}\frac{\phi(n_i+k-1)}{\phi(n_i+k)} = \alpha_{k-1}.
\]
This is rather stronger than the statement that \eqref{compare.Schinzel} occurs infinitely often. However, our corollary (in addition to being a specialization of a theorem on more general, not necessarily monic linear polynomials) has the advantage that it provides a set of solutions to \eqref{compare.Schinzel} of positive lower density.

We can draw another corollary from Theorem \ref{general.theorem}. Given integers $a,b,c,d$ with $a,b>0$ and $ad\ne bc$, we can take $a_1=\dots=a_k=a$ and $b,b_1,\dots,b_k$ in arithmetic progression with common difference $a$, and similarly $a_{k+1}=\dots=a_{2k}=c$ and $d,b_{k+1},\dots,b_{2k}$ in arithmetic progression with common difference $c$ in Theorem \ref{general.theorem}, to obtain a result on consecutive solutions to the fixed inequality $\phi(an+b) > \phi(cn+d)$:

\begin{cor}
Let $a$ and $b$ be positive integers and $c$ and $d$ any integers such that $ad\ne bc$. Then for any positive integer $k$, the set of positive integers $n$ for which
\begin{multline}
\phi\big(a(n+1)+b\big) > \phi\big(c(n+1)+d\big),\quad \phi\big(a(n+2)+b\big) > \phi\big(c(n+2)+d\big), \\
\dots,\quad \phi\big(a(n+k)+b\big) > \phi\big(c(n+k)+d\big)
\label{consecutive.ineqs}
\end{multline}
has positive lower density.
\label{consecutive.cor}
\end{cor}

This corollary extends a theorem of Newman \cite{newman} on individual solutions to arbitrarily long strings of consecutive solutions, and also quantifies the set of solutions to have positive lower density. Presumably the set of solutions to any system of inequalities of the form \eqref{general.ineqs}, \eqref{compare.Schinzel}, or \eqref{consecutive.cor} actually possesses a well-defined density, but since we do not have a proof that this property holds, we simply show that any such set of solutions has positive lower density, that is, contains a subset of positive density.

As it happens, the sum-of-divisors function $\sigma(n) = \sum_{d\mid n} d$ is closely enough related to $\phi(n)$ that the truth of the results mentioned above is enough to imply the analogous results for $\sigma$:

\begin{cor}
Theorem \ref{general.theorem} and Corollaries \ref{permutation.cor} and \ref{consecutive.cor} hold with the function $\phi$ replaced by $\sigma$ everywhere.
\label{sigma.cor}
\end{cor}

\begin{proof}
Given that
\[
\frac{\phi(m)}m \frac{\sigma(m)}m = \prod_{p\mid m} \bigg( 1-\frac1p \bigg) \prod_{p^\alpha\|m} \bigg( 1+\frac1p+\dots+\frac1{p^\alpha} \bigg) = \prod_{p^\alpha\|m} \bigg( 1 - \frac1{p^{\alpha+1}} \bigg),
\]
we see immediately that
\[
1 \ge \frac{\phi(m)\sigma(m)}{m^2} > \prod_p \bigg( 1-\frac1{p^2} \bigg) = \frac1{\zeta(2)} = \frac6{\pi^2},
\]
or in other words
\[
\frac{m^2}{\phi(m)} \ge \sigma(m) > \frac{6m^2}{\pi^2\phi(m)}.
\]
If we know that $\phi(a_jn+b_j)/\phi(a_{j+1}n+b_{j+1}) > C$, it therefore follows that for any fixed $\ep>0$, we have
\begin{align*}
\frac{\sigma(a_{j+1}n+b_{j+1})}{\sigma(a_jn+b_j)} &> \frac{6(a_{j+1}n+b_{j+1})^2/\pi^2\phi(a_{j+1}n+b_{j+1})}{(a_jn+b_j)^2/\phi(a_jn+b_j)} \\
&= \frac6{\pi^2} \frac{\phi(a_jn+b_j)}{\phi(a_{j+1}n+b_{j+1})} \frac{(a_{j+1}n+b_{j+1})^2}{(a_jn+b_j)^2} > \frac{6C}{\pi^2}\bigg( \bigg( \frac{a_{j+1}}{a_j} \bigg)^2 - \ep \bigg)
\end{align*}
when $n$ is sufficiently large in terms of $\ep$. Therefore if we choose any positive constant $C'$ with
\[
C' < \frac{6C}{\pi^2} \bigg( \frac{\min\{a_1,\dots,a_k\}}{\max\{a_1,\dots,a_k\}} \bigg)^2,
\]
we see that all sufficiently large solutions $n$ of the inequalities \eqref{general.ineqs} also satisfy
\[
\frac{\sigma(a_kn+b_k)}{\sigma(a_{k-1}n+b_{k-1})} > C', \quad\dots, \quad \frac{\sigma(a_3n+b_3)}{\sigma(a_2n+b_2)} > C', \quad \frac{\sigma(a_2n+b_2)}{\sigma(a_1n+b_1)} > C'.
\]
Since $C$ can be taken arbitrarily large in Theorem \ref{general.theorem}, we may take $C'$ as large as we want here. This proves the analogue of Theorem \ref{general.theorem} for the function $\sigma$ (upon reversing the orders of the $a_j$ and $b_j$). Given this analogue for $\sigma$, the analogues of Corollaries \ref{permutation.cor} and \ref{consecutive.cor} for $\sigma$ follow in exactly the same way as for $\phi$.
\end{proof}

Now that we have finished showing that the three corollaries above follow from Theorem \ref{general.theorem}, we lay down our plan of attack for the proof of that theorem. The idea is to restrict to an arithmetic progression carefully constructed to be extremely biased in favor of the inequalities \eqref{general.ineqs}. Once we have such an arithmetic progression, we then show that a positive proportion of numbers in that progression actually do satisfy the inequalities as expected.

We begin by proving two lemmas. The first lemma is simply a bit of bookkeeping that it will be convenient to have in hand when we attack Theorem \ref{general.theorem} in earnest.

\begin{lemma}
Let $a_1,\dots,a_k,b_1,\dots,b_k$ be positive integers with $a_ib_j\ne a_jb_i$ for all $1\le i<j\le k$, and let $A$ be a nonzero integer that is a multiple of each $b_j$ $(1\le j\le k)$ and also a multiple of each $a_ib_j-a_jb_i$ $(1\le i<j\le k)$. Let $Q_1,\dots,Q_k$ be pairwise relatively prime positive integers, each of which is also relatively prime to $A$. Finally, let $r$ be an integer such that $Q_j \mid (a_jAr+b_j)$ for each $1\le j\le k$. Then for each $1\le j\le k$, the greatest common divisor of $a_jAQ_1Q_2\cdots Q_k$ and $a_jAr+b_j$ is exactly $Q_jb_j$.
\label{gcd.lemma}
\end{lemma}

\begin{proof}
Suppose first that $1\le i,j\le k$ with $i\ne j$ and that $p$ is a prime factor of $Q_i$. Then $a_iAr+b_i\equiv 0\mod p$ by the hypothesis on $r$. If it were also the case that $a_jAr+b_j\equiv 0\mod p$, then we would have
\[
a_jb_i - a_ib_j = a_j(a_iAr+b_i) - a_i(a_jAr+b_j) \equiv 0\mod p,
\]
and so $p$ divides $A$, contradicting the coprimality of $A$ and $Q_i$. Therefore $p$ does not divide $a_jAr+b_j$. Since this is true of every prime dividing the $Q_i$ with $i\ne j$, we may ignore the factor $Q_1\dots Q_{j-1}Q_{j+1}\dots Q_k$ when computing greatest common divisors with $a_jAr+b_j$. In other words, the greatest common divisor of $a_jAQ_1Q_2\cdots Q_k$ and $a_jAr+b_j$ is the same as the greatest common divisor of $a_jAQ_j$ and $a_jAr+b_j$.

It is then immediate that this greatest common divisor divides the linear combination $Q_j(a_jAr+b_j) - r(a_jAQ_j) = Q_jb_j$. On the other hand, $Q_j$ divides $a_jAQ_j$ obviously and divides $a_jAr+b_j$ by hypothesis; also, $b_j$ divides both $a_jAQ_j$ and $a_jAr+b_j$ as well, since $A$ is a multiple of $b_j$. Moreover, $Q_j$ and $b_j$ are relatively prime, since $b_j$ divides $A$ which is relatively prime to $Q_j$ by assumption. Therefore the greatest common divisor in question is also a multiple of $Q_jb_j$, which establishes the lemma.
\end{proof}

The second of the two lemmas is somewhat more interesting: it asserts that given a collection of suitable linear functions, there are many values of the input variable for which all of the function values have relatively large $\phi$-value. The method of proof is based upon a suggestion of Pomerance. We have chosen to use simple constants in the lemma rather than increase the technical demands by attempting to optimize them.

\begin{lemma}
Let $k$ be a positive integer, and let $c_1,\dots,c_k,d_1,\dots,d_k$ be positive integers with $c_j$ and $d_j$ relatively prime for each $1\le j\le k$. Then the set of positive integers $m$ for which the inequalities
\[
\frac{\phi(c_1m+d_1)}{c_1m+d_1} > e^{-k},\quad \frac{\phi(c_2m+d_2)}{c_2m+d_2} > e^{-k},\quad \dots,\quad \frac{\phi(c_km+d_k)}{c_km+d_k} > e^{-k}
\]
are all satisfied has lower density at least $\frac3{10}$.
\label{wise.lemma}
\end{lemma}

\begin{proof}
Suppose the statement of the lemma were false. Then on a set of positive integers $m$ of upper density greater than $\frac7{10}$, at least one of the factors $\phi(c_jm+d_j)/(c_jm+d_j)$ is at most $e^{-k}$, that is,
\[
\log \bigg( \frac{\phi(c_jm+d_j)}{c_jm+d_j} \bigg)^{\!-1} = \sum_{p\mid (c_jm+d_j)} \log\bigg( 1-\frac1p\bigg)^{\!-1} \ge k
\]
for at least one index $1\le j\le k$. In particular, for any such integer $m$ we have
\[
\sum_{j=1}^k \sum_{p\mid (c_jm+d_j)} \log\bigg( 1-\frac1p\bigg)^{\!-1} \ge k
\]
since all the terms are nonnegative. Since the set of such integers $m$ has upper density greater than $\frac7{10}$, there are arbitrarily large values of $x$ up to which there are at least $\frac7{10}x$ such integers $m$, so that
\begin{equation}
\sum_{m\le x} \sum_{j=1}^k \sum_{p\mid (c_jm+d_j)} \log\bigg( 1-\frac1p\bigg)^{\!-1} \ge \frac{7kx}{10}.
\label{lower.bound.contr}
\end{equation}

On the other hand, using the fact that $\log(1-x)^{-1} \le x\log 4$ for $0<x\le\frac12$, we have
\[
\sum_{m\le x} \sum_{j=1}^k \sum_{p\mid (c_jm+d_j)} \log\bigg( 1-\frac1p\bigg)^{\!-1} \le \sum_{j=1}^k \sum_{p \le xM} \sum_{\substack{m\le x \\ p\mid (c_jm+d_j)}} \frac{\log 4}p,
\]
where we have defined $M = 2\max\{c_1,\dots,c_k,d_1,\dots,d_k\}$; note in particular that $M$ does not depend on $x$. Since $c_j$ and $d_j$ are relatively prime, no prime $p$ that divides $c_j$ ever divides $c_jm+d_j$. Therefore at most one value of $m$ in every set of $p$ consecutive values makes $c_jm+d_j$ a multiple of $p$, and so we have
\begin{align*}
\sum_{j=1}^k \sum_{p\le xM} \sum_{\substack{m\le x \\ p\mid (c_jm+d_j)}} \frac{\log 4}p &\le \sum_{j=1}^k \sum_{p \le xM} \frac{\log 4}p \bigg( \frac xp + 1 \bigg) \\
&\le (k \log 4) \bigg( x \sum_p \frac1{p^2} + \sum_{2\le \ell\le xM} \frac1\ell \bigg) \\
&\le (k \log 4) \bigg( \frac x2 + \log xM \bigg),
\end{align*}
where we have used the fact that $\sum_p 1/p^2 < 1/2$. We conclude that
\begin{equation}
\sum_{m\le x} \sum_{j=1}^k \sum_{p\mid (c_jm+d_j)} \log\bigg( 1-\frac1p\bigg)^{\!-1} \le(k\log 2)x + (k\log 4)\log xM.
\label{upper.bound.contr}
\end{equation}

Combining the inequalities \eqref{lower.bound.contr} and \eqref{upper.bound.contr} yields
\[
\frac{7kx}{10} \le(k\log 2)x + (k\log 4)\log xM
\]
for arbitrarily large values of $x$. However, since $\log 2 < \frac7{10}$, this results in a contradiction when $x$ is sufficiently large. This contradiction establishes the lemma.
\end{proof}

We are now prepared to establish the main theorem of this paper.

\begin{proof}[Proof of Theorem \ref{general.theorem}]
By replacing $n$ with $n+n_0$ for sufficiently large $n_0$, we may assume that in fact all of the integers $b_1,\dots,b_k$ are positive; this changes the number of solutions by only a finite amount.

Define
\[
\nu = \min\bigg\{ \frac{\phi(b_1)}{b_1}, \dots, \frac{\phi(b_k)}{b_k} \bigg\} \quad\text{and}\quad \rho = \frac{\min\{a_1,\dots,a_k\}}{2\max\{a_1,\dots,a_k\}},
\]
and choose any number $\ep$ satisfying
\begin{equation}
0 < \ep < \frac{e^{-k}\nu\rho}{2C}.  \label{epsilon}
\end{equation}
Choose a positive integer $A$ such that all of the following numbers divide $A$: each of $a_1,\dots,a_k$, each of $b_1,\dots,b_k$, and every $a_ib_j-a_jb_i$ for $1\le i<j\le k$. Define $\A$ to be the finite set of primes that divide $A$. Now, choose sets of primes $\P_1,\dots,\P_k$, pairwise disjoint and also disjoint from $\A$, that satisfy the inequalities
\begin{equation}
\tfrac12\ep^{j-1} < \prod_{p\in\P_j} \bigg( 1-\frac1p \bigg) \le \ep^{j-1}
\label{P.inequalities}
\end{equation}
for each $1\le j\le k$; this is possible since the product $\prod \big(1-\frac1p\big)$ diverges to 0 and the individual factors $1-\frac1p$ are always between $\frac12$ and 1. Notice that we might well take $\P_1$ to be the empty set, although the precise choice of the $\P_j$ plays little role in the arguments to follow.

Define $Q_j = \prod_{p\in\P_j} p$ for each $1\le j\le k$, and define $Q = Q_1Q_2\dots Q_k$. Define the positive integer $r$ to be a solution to the system of congruences
\begin{equation}
a_jAr + b_j \equiv 0\mod{Q_j} \quad (1\le j\le k).
\label{CRT.congruences}
\end{equation}
Each individual congruence specifies a unique $r\mod{Q_j}$: the coefficient $a_jA$ is composed entirely of primes in the set $\A$ and so cannot have any factor in common with $Q_j$, since the prime factors of $Q_j$ lie in the set $\P_j$ which is disjoint from $\A$. By the Chinese Remainder Theorem, an integer $r$ exists (and is unique modulo $Q$) that satisfies all the congruences \eqref{CRT.congruences} simultaneously.

We now consider the $k$ linear polynomials $a_jx+b_j$ evaluated on the arithmetic progression $x = A(Qm+r)$, that is, we consider the polynomials $a_jAQm + (a_jAr+b_j)$ with $m$ as the variable. Recalling that $Q=Q_1Q_2\dots Q_k$, we see that the numbers $a_j,b_j,Q_j,A,r$ satisfy the hypotheses of Lemma \ref{gcd.lemma}, and so we have $(a_jAQ,a_jAr+b_j) = Q_jb_j$ for every $1\le j\le k$. We may therefore write $a_j(A(Qm+r))+b_j$ as $Q_jb_j(c_jm+d_j)$, where $c_j$ and $d_j$ are relatively prime integers.

Using the elementary inequality $\phi(st) \ge \phi(s)\phi(t)$, we can derive the lower bound
\begin{align*}
\frac{\phi(a_jA(Qm+r)+b_j)}{a_jA(Qm+r)+b_j} &= \frac{\phi(Q_jb_j(c_jm+d_j))}{Q_jb_j(c_jm+d_j)} \\
&\ge \frac{\phi(Q_j)}{Q_j} \frac{\phi(b_j)}{b_j} \frac{\phi(c_jm+d_j)}{c_jm+d_j} \ge \tfrac12\ep^{j-1} \cdot \nu \cdot \frac{\phi(c_jm+d_j)}{c_jm+d_j}.
\end{align*}
Define
\[
\M = \bigg\{ m\ge1\colon \frac{\phi(c_1m+d_1)}{c_1m+d_1} > e^{-k},\quad \frac{\phi(c_2m+d_2)}{c_2m+d_2} > e^{-k},\quad \dots,\quad \frac{\phi(c_km+d_k)}{c_km+d_k} > e^{-k} \bigg\}.
\]
By Lemma \ref{wise.lemma}, the set $\M$ has lower density at least $\frac3{10}$. For $m\in\M$, we have
\begin{equation}
\frac{\phi(a_jA(Qm+r)+b_j)}{a_jA(Qm+r)+b_j} > \tfrac12\ep^{j-1} \cdot \nu \cdot e^{-k} > \frac{C\ep^j}\rho
\label{lower.bound}
\end{equation}
by the choice \eqref{epsilon} of $\ep$.

A complementary upper bound is even easier to derive, using the general inequality $\phi(st) \le \phi(s)t$:
\begin{equation}
\frac{\phi(a_{j+1}A(Qm+r)+b_{j+1})}{a_{j+1}A(Qm+r)+b_{j+1}} = \frac{\phi(Q_{j+1}b_{j+1}(c_{j+1}m+d_{j+1}))}{Q_{j+1}b_{j+1}(c_{j+1}m+d_{j+1})} \le \frac{\phi(Q_{j+1})}{Q_{j+1}} \le \ep^j.
\label{upper.bound}
\end{equation}

Now let $\N = \{A(Qm+r)\colon m\in\M\}$, which has lower density at least $3/10AQ$. For any $n\in\N$, the inequality \eqref{lower.bound} tells us that $\phi(a_jn+b_j) > C\ep^j(a_jn+b_j)/\rho$, while the inequality \eqref{upper.bound} tells us that $\phi(a_{j+1}n+b_{j+1}) \le \ep^j(a_{j+1}n+b_{j+1})$. Therefore
\begin{align*}
\frac{\phi(a_jn+b_j)}{\phi(a_{j+1}n+b_{j+1})} &> \frac{C\ep^j(a_jn+b_j)/\rho}{\ep^j(a_{j+1}n+b_{j+1})} \\
&= C \frac{a_j/a_{j+1}}\rho \frac{n+b_j/a_j}{n+b_{j+1}/a_{j+1}} \ge 2C \frac{n+b_j/a_j}{n+b_{j+1}/a_{j+1}}
\end{align*}
by the definition of $\rho$. The last fraction tends to 1 as $n$ tends to infinity, so we conclude that
\[
\frac{\phi(a_jn+b_j)}{\phi(a_{j+1}n+b_{j+1})} > C
\]
for all $n\in\N$ that are sufficiently large in terms of the $a_j$ and $b_j$. Since this set of $n$ has positive lower density, the theorem is proved.
\end{proof}

\bigskip{\small
{\it Acknowledgements.} The author thanks Michael Filaseta for asking the thought-provoking question that led to the investigation of these results, as well as Carl Pomerance for a wise suggestion that led to the successful use of Lemma \ref{wise.lemma}, and also expresses his appreciation to the Centre de recherches math\'ematiques under whose hospitality he completed this manuscript. The author was supported in part by a Natural Sciences and Engineering Research Council grant.}

\bibliographystyle{amsplain}
\bibliography{phinequalities}

\end{document}